\documentclass[a4paper,12pt,reqno]{amsart}
\usepackage{latexsym}
\usepackage{amsmath,amsthm,amssymb,amscd}
\usepackage{fullpage}
\usepackage{xcolor}
\usepackage{parskip}
\usepackage{mathtools}
\mathtoolsset{showonlyrefs}

\usepackage[activate={true,nocompatibility},final,tracking=true,kerning=true,spacing=true]{microtype}
\usepackage[pdfencoding=unicode,pdftex,bookmarks=true,bookmarksnumbered]{hyperref}
\definecolor{citecolour}{rgb}{0.11, 0.22, 0.73}
\definecolor{linkcolour}{rgb}{0.03, 0.47, 0.18}
\hypersetup{colorlinks,breaklinks,
		linkcolor=linkcolour,citecolor=citecolour,
		filecolor=linkcolour, urlcolor=linkcolour}
\usepackage{version}

\theoremstyle{plain}
\newtheorem{theorem}{Theorem}[section]

\newtheorem{lemma}[theorem]{Lemma}

\newtheorem{corollary}[theorem]{Corollary}

\theoremstyle{remark}

\theoremstyle{definition}

\newtheorem*{notation*}{Notation}

\numberwithin{equation}{section}

\DeclareMathOperator{\M}{M}
\DeclareMathOperator{\J}{J}
\DeclareMathOperator{\HS}{HS}
\DeclareMathOperator{\Suz}{Suz}
\DeclareMathOperator{\McL}{McL}
\DeclareMathOperator{\Ru}{Ru}
\DeclareMathOperator{\He}{He}
\DeclareMathOperator{\Ly}{Ly}
\DeclareMathOperator{\ON}{O'N}
\DeclareMathOperator{\Co}{Co}
\DeclareMathOperator{\Fi}{Fi}
\DeclareMathOperator{\HN}{HN}
\DeclareMathOperator{\Th}{Th}

\DeclareMathOperator{\baby}{\mathbb{B}}
\DeclareMathOperator{\monster}{\mathbb{M}}
\DeclareMathOperator{\psl}{PSL}
\DeclareMathOperator{\pgl}{PGL}
\DeclareMathOperator{\psu}{PSU}
\DeclareMathOperator{\agam}{A\Gamma L}
\DeclareMathOperator{\agl}{AGL}
\DeclareMathOperator{\psp}{PSp}
\DeclareMathOperator{\Sp}{Sp}

\DeclareMathOperator{\Zent}{\mathbf{Z}}

\DeclareMathOperator{\Aut}{\mathbf{Aut}}
\DeclareMathOperator{\GL}{GL}
\DeclareMathOperator{\cent}{\mathbf{C}}
\DeclareMathOperator{\fitt}{\mathbf{F}}
\DeclareMathOperator{\fratt}{\mathbf{\Phi}}
\DeclareMathOperator{\Soc}{\mathbf{Soc}}
\renewcommand{\leq}{\leqslant}
\renewcommand{\geq}{\geqslant}
\def\GAP{{\sf GAP}}

\newenvironment{proofof}{{\bf {Proof.} }}{\hfill $\blacksquare$ \\}
\newenvironment{proofofA}{{\bf {Proof of Theorem 2.1.} }}{\hfill $\blacksquare$ \\}

\begin{document}
\title{A note on abelian subgroups of maximal order}
\date{\today}
\author{Stefanos Aivazidis}
\email{stefanosaivazidis@gmail.com}
\author{Robert M. Guralnick}
\address{Department of Mathematics, University of Southern California,
Los Angeles, CA 90089-2532, USA}
\email{guralnic@usc.edu}
\thanks{The second author was partially supported by the NSF
grant DMS-1600056.}
\maketitle
\section{Introduction}
In a recent paper of the first author and I. M. Isaacs it was shown that 
if $m=m(G)$ is the maximal order of an abelian subgroup of the finite group $G$, 
then $|G|$ divides $m!$ (\cite[Thm. 5.2]{manl}).
The purpose of this brief note is to improve on the $m!$ bound (see Theorem~\ref{thm:A} below). 
We shall then take up the task of determining when the (implicit) inequality of our theorem becomes an equality.
Despite, perhaps, first appearances this determination is not trivial. 
To accomplish it we shall establish a result (Theorem~\ref{thm:B}) of independent interest and we shall then see that 
Theorems \ref{thm:A} and \ref{thm:B} combine to further strengthen Theorem~\ref{thm:A} (see Theorem~\ref{thm:g/h}).
\section{Theorems \ref{thm:A} and \ref{thm:B}}
\begin{theorem}\label{thm:A}
For each positive integer $n$, define $g(n)$ to be the product of all prime powers at most $n$. 
Let $G$ be a finite group and suppose that $|A| \leq m$ for every abelian subgroup $A$ of $G$. Then $|G|$ divides $g(m)$.
\end{theorem}
For the proof we shall require an auxiliary lemma and we present that first. 
We remark that the lemma is a direct consequence of a classic result due to Burnside \cite{burnside} which asserts that 

\textit{if $P$ has order $p^n$,  $\Zent(P)$ has order $p^c$
and $p^v$ is the maximal order of a normal abelian subgroup of $P$, then
$n - v \leq 1/2(v-c)(v + c-1)$}

upon taking $c=1$. 
Nevertheless, we choose to present our own proof.

\begin{lemma}\label{lem}
If $G$ is a $p$-group of order $p^k$ and $A$ is a maximal abelian normal subgroup of $G$, where $|A|=p^s$, then $k \leq s(s+1)/2$.
\end{lemma}
\begin{proofof}
We assume familiarity with the following well-known theorem of Hall:

\emph{If $G$ is finite and can be generated by $r$ elements then}
\begin{equation}
\left\lvert \Aut(G)\right\rvert \,\,\, \mbox{\emph{divides}}\,\,\, \left\lvert \fratt(G)\right\rvert^r \left\lvert \Aut(G/\fratt(G)) \right\rvert \, .
\end{equation}
As a consequence of this theorem of Hall and Burnside's Basis Theorem, 
if $G$ is a $p$-group of order $p^k$ and $p^d = |G : \fratt(G)|$ then
\[
\lvert \Aut(G) \rvert \, \Big| \,\, p^{d(k-d)} \left\lvert \GL_d(p)\right\rvert  = p^{d(k-d)} \prod_{i=0}^{d-1}(p^d-p^i)\, .
\]
Thus a Sylow $p$-subgroup of $\Aut(G)$ has order at most 
\begin{equation}\label{bound}
p^{d(k-d)+d(d-1)/2} \leq p^{k(k-1)/2}\,.
\end{equation}
To prove the claim, it suffices to prove that $|G:A| \leq p^{s(s-1)/2}$. 
Since $G$ is nilpotent, every maximal abelian normal subgroup of $G$ is self-centralising (cf. \cite[Lemma 4.16]{isaacs}). 
In particular, the $N/C$ Theorem shows that $|G:A|$ divides the $p$-part of $\lvert\Aut(A)\rvert$, 
which is at most $p^{s(s-1)/2}$ by \eqref{bound}.
\end{proofof}

\begin{proofofA}
It will suffice to prove that $|P|$ divides $g(m)$ for each Sylow subgroup $P$ of $G$. 
Fix $P$ and let $A$ be a maximal abelian normal subgroup of $P$. 
Suppose $|A|=p^s$. Lemma~\ref{lem} shows that $|P| \leq p^{s(s+1)/2}$, thus $|P|$ divides $p^{s(s+1)/2}$. 
Note that
\begin{equation}
p^{s(s+1)/2} = (g(|A|))_p \,\, \big| \,\, g(|A|) \,\, \big| \,\, g(m) \, ,
\end{equation}
where $g(|A|)$ divides $g(m)$ since $|A| \leq m$. 
Thus $|P|$ divides $g(m)$, as wanted.
\end{proofofA}

\begin{theorem}\label{thm:B}
Let $G$ be a finite group and suppose that $m$ is the maximal order
of an abelian subgroup of $G$. Suppose that $|G|$ is divisible by a prime 
$p> m/2$. Then one of the following holds:
\begin{enumerate}
\item  $G=PC$ where $|P|=p$ and $C$ acts faithfully on $P$ (and so is
non-trivial cyclic of order dividing $p-1$); 
\item   $G \cong S_3$ with $p=2$ or $3$;
\item   $G$ contains $\fitt( \agam_1 (2^a))$ where  
    $p=2^a-1$ is a Mersenne prime; or 
\item   $G$ is almost simple.
\end{enumerate}
\end{theorem}

\begin{proofof}  Let $P$ be a Sylow $p$-subgroup of $G$.  Note that $P$ has
order $p$ and is self-centralising. 
Let $N$ be a minimal normal subgroup of $G$.   If $p$ does not divide $|N|$, then
$P$ acts fixed point freely on $N$ and so $N$ is nilpotent (by Thompson's thesis that
Frobenius kernels are nilpotent) and $|N| > p$.   Thus, 
$N$ is self-centralising and so $N=\fitt(G)=\fitt^*(G)$.   Suppose that $N$ is elementary
abelian of order $r^a$ for some prime $r \ne p$. If $a=1$ then $p$
divides $r-1$ and $p > r/2$ and so $p=r-1$, whence $r =3$ and $p=r-1=2$. Therefore, 
$G \cong S_3$ and we are in case (2).

Now assume that $a \geq 2$. Since $p > r^a/2$, this implies
that $p=r^a-1$, whence $r=2$ and $p$ is a Mersenne prime.  Thus,  $G=NH$
where $H$ is a subgroup of $\GL_a(2)$ containing a Singer cycle (i.e. an element
of order $2^a-1$).  Note that since $2^a-1$ is prime, this implies
that $a$ is prime.  By \cite{kantor},  this implies, in turn, that either $H=\GL_a(2)$ or that 
$G \leq  \agam_1(2^a)$.    The last case is allowed. Notice, here, that if $G \leq  \agam_1(2^a)$ then
also $G \geq N= \fitt( \agam_1(2^a))$.   If $a =2$,
   then $ S_4 \cong \agam_1(2^2)$ and $G=A_4$ or $S_4$ 
and the result follows.
 If $a=3$, then either $G \leq \agam_1(2^3)$ , or $G=\agl_3(2)$.
We discount the latter possibility by noting that $\lvert \agl_3(2) \rvert = 2^6.3.7$, while $m(\agl_3(2))=16$.

So assume that $a > 3$.    Note that $\GL_a(2)$ contains an abelian subgroup of 
order $2^{(a^2-1)/4}$.   This is greater than $2p$ for $a \geq 5$ and so only (3)
occurs in this case.   

Now suppose that $p$ divides $|N|$.  If $N$ is a $p$-group, then (1) holds.
Otherwise $N$ must be a non-abelian simple group and $\cent_G(N)=1$.  Thus,
$G$ is almost simple with socle $N$.
\end{proofof}

Now we check that the groups in cases (1), (2), and (3) have a  unique
prime divisor of order greater than $m/2$.   Thus we have:

\begin{corollary} \label{2prime}   Let $G$ be a finite group and suppose that $m$ 
is the maximal order
of an abelian subgroup of $G$. Suppose that $|G|$ is divisible by primes $p, q$
with $p> q > m/2$.   Then either $m=p=3$ and  $G \cong S_3$ or  $G$ is almost simple.
\end{corollary}

We shall later see (in Theorem~\ref{thm:C} and Corollary~\ref{cor:first}) that, in fact, 
there are no genuine almost simple groups satisfying the conclusion of Corollary~\ref{2prime} 
and so the only such groups are $S_3$ and certain simple non-abelian groups.

\section{Large primes of almost simple groups}
We now attempt to classify the simple groups $G$ whose order is divisible by
one or two ``large" primes
$p$, $q$ (i.e. such that $m/2 < p < q$), and we begin with the sporadic simple groups. For those
sporadic groups whose order is relatively small, we can compute $m(G)$ explicitly using {\GAP}~\cite{GAP}. 
When that is either time-consuming or just impossible without resorting to theoretical arguments, 
we are content to provide a lower bound for $m(G)$ by using Lemma~\ref{lem} in the form

\textit{if $w$ is the smallest positive integer such that $1+\ldots +w \geq a$, 
then a group of order $r^a$, $r$ a prime, has an abelian (normal) subgroup of order at least $r^w$}

and considering the maximum $w$ over all primes involved in the factorisation of $|G|$. 
This only fails to provide a sufficient lower bound for $m(G)$ when $G \cong \Ly$. For that group
we consult \cite[p. 174]{atlas} and find that it has a maximal subgroup of shape $3^5{:}(2 \times \M_{11})$,
whence $m(\Ly) \geq 3^5$.

Our findings, which are summarised in Table~\ref{tab:sporadic}, 
indicate that the only sporadic simple groups with exactly two large prime divisors are $\J_1$ and $\J_3$, 
whereas $\M_{11}$, $\M_{12}$, $\M_{22}$, and $\M_{23}$ each have a single large prime. 
An extra word on how Table~\ref{tab:sporadic} is set-up: $p$ is simply the largest prime divisor of $|G|$. 
It is clear that $|G|$ has no large prime divisors if and only if $m(G) \geq 2p$.

\begin{table}[ht!]
\centering
\begin{tabular}{c c c c}
\hline
$G$ &  $|G|$ & $m(G)$ & $2p$ \\
\hline
$\M_{11}$ 	& $2^4.3^2.5.11$ 										& $=11$ 		& $22$\\
$\M_{12}$ 	& $2^6.3^2.5.11$ 										& $=2^4$ 	& $22$\\
$\M_{22}$ 	& $2^7.3^2.5.7.11$ 										& $=2^4$ 	& $22$\\
$\M_{23}$ 	& $2^7.3^2.5.7.11.23$ 									& $=23$ 		& $46$\\
$\M_{24}$ 	& $2^{10}.3^3.5.7.11.23$ 								& $=2^6$ 	& $46$\\ 
$\J_1$ 		& $2^3.3.5.7.11.19$ 									& $=19$ 		& $38$\\
$\J_2$ 		& $2^7.3^3.5^2.7$ 										& $=5^2$ 	& $14$\\
$\J_3$ 		& $2^7.3^5.5.17.19$ 									& $=3^3$ 	& $38$\\
$\J_4$ 		& $2^{21}.3^3.5.7.11^3.23.29.31.37.43$ 					& $\geq 11^2$	& $86$\\
$\HS$ 		& $2^9.3^2.5^3.7.11$ 									& $=2^6$ 	& $22$\\
$\McL$ 		& $2^7.3^6.5^3.7.11$ 									& $=3^4$ 	& $22$\\
$\He$ 		& $2^{10}.3^3.5^2.7^3.17$ 								& $=2^6$		& $34$\\
$\Ru$ 		& $2^{14}.3^3.5^3.7.13.29$ 								& $= 2^6$ 	& $58$\\
$\Suz$ 		& $2^{13}.3^7.5^2.7.11.13$ 								& $= 3^5$ 	& $26$\\
$\ON$ 		& $2^9.3^4.5.7^3.11.19.31$ 								& $= 3^4$ 	& $62$\\
$\Co_3$ 		& $2^{10}.3^7.5^3.7.11.23$ 								& $= 3^5$ 	& $46$\\
$\Co_2$ 		& $2^{18}.3^6.5^3.7.11.23$ 								& $\geq 2^6$ 	& $46$\\
$\Co_1$ 		& $2^{21}.3^9.5^4.7^2.11.13.23$ 							& $\geq 3^4$ 	& $46$\\
$\Fi_{22}$ 	& $2^{17}.3^9.5^2.7.11.13$ 								& $\geq 3^4$ 	& $26$\\
$\Fi_{23}$ 	& $2^{18}.3^{13}.5^2.7.11.13.17.23$ 						& $\geq 3^5$ 	& $46$\\
$\Fi_{24}'$ 	& $2^{21}.3^{16}.5^2.7^3.11.13.17.23.29$ 					& $\geq 3^6$ 	& $58$\\
$\HN$ 		& $2^{14}.3^6.5^6.7.11.19$ 								& $\geq 5^3$ 	& $38$\\
$\Ly$ 		& $2^8.3^7.5^6.7.11.31.37.67$ 							& $\geq 3^5$ 	& $134$\\
$\Th$ 		& $2^{15}.3^{10}.5^3.7^2.13.19.31$ 						& $\geq 3^4$ 	& $62$\\
$\baby$ 		& $2^{41}.3^{13}.5^6.7^2.11.13.17.19.23.31.47$ 				& $\geq 2^9$ 	& $94$\\
$\monster$ 	& $2^{46}.3^{20}.5^9.7^6.11^2.13^3.17.19.23.29.31.41.47.59.71$ 	& $\geq 2^{10}$ & $142$\\
\hline
\end{tabular}
\caption{The sporadic simple groups $G$}
\label{tab:sporadic}
\end{table}

From this point onward our references regarding estimates for the maximal size of an abelian subgroup will be \cite{vdovin} and \cite[Table 3.3.1]{gls3}.
\subsection{Alternating case}
We consider alternating groups first. Let $n=3k+r$, $0 \leq r \leq 2$. Clearly, $A_n$ has an abelian subgroup of order $3^k \geq 3^{\frac{n-2}{3}}$. If we denote by $p$ the largest prime divisor of $|A_n|$, then $p\leq n$. A straightforward induction on $n$ shows that 
$3^{\frac{n-2}{3}} > 2n$ for all $n\geq 11$, thus $m(A_n) > 2p$ for all $n$ in that range. For $5 \leq n \leq 10$, we find $2p=10$ for $A_5$ and $A_6$, while $2p=14$ for $n>6$. On the other hand, $m(A_5)=5$, $m(A_7)=12$, and $m(A_8)=16$ whence $A_5$ is the only alternating group with two large prime divisors, while $A_6$ and $A_7$ have a single large prime divisor.
\subsection{Linear case}
Next in line are the classical groups and we begin with the linear case. 
For $n \geq 4$ we have $m(\psl_n(q)) = q^{\lfloor \frac{n^2}{4}\rfloor}$. We recall that
\begin{equation}
\lvert\psl_n(q)\rvert = \frac{1}{\gcd(n,q-1)}q^{n(n-1)/2}\prod_{i=2}^{n}(q^i-1)\,.
\end{equation}
It is easy to see that the largest prime divisor of $\lvert\psl_n(q)\rvert$
is at most $\frac{q^n-1}{q-1}$. For $n \geq 5$, we have $\lfloor \frac{n^2}{4}\rfloor \geq n+1$, thus 
\begin{equation}
m(\psl_n(q)) \geq \frac{q^{n+1}}{2} \geq q^n > \frac{q^n-1}{q-1}\,.
\end{equation}
Therefore, it suffices to consider the cases arising from $n \in \{2,3,4\}$. 
The case $n=4$ can be dispensed with quickly. 
For $1+q+q^2+q^3= (q^2+1)(q+1)$, thus the largest prime divisor of $\lvert\psl_4(q)\rvert$ cannot exceed 
$1+q+q^2 = \frac{q^3-1}{q-1}$; but $\frac{q^4}{2} > 1+q+q^2$ for all prime powers $q$.
For $n=3$ we have
\begin{equation}
\lvert\psl_3(q)\rvert = \frac{1}{\gcd(3,q-1)}q^{3}(q-1)^2(q+1)(q^2+q+1)\,.    
\end{equation}
Moreover, $m(\psl_3(q)) = q^2+q+1$ if $\gcd(3,q-1)=1$ and $q^2$ otherwise. 

Suppose first that $\gcd(3,q-1)=1$. 
If $q^2+q+1$ is a prime, then $m=q^2+q+1$ is the largest prime divisor of the order. 
In that case, the second largest prime divisor is at most $q+1 < \frac{q^2+q+1}{2} = \frac{m}{2}$. 
If, on the other hand, $q^2+q+1$ is not a prime, then the largest prime divisor of the order 
is at most $\max\left\{\frac{q^2+q+1}{3}, q+1\right\} < \frac{m}{2}$.   So there are never two
large prime divisors and there is a single large prime divisor precisely when 
$$q^2+q+1 = p =m.$$

Now suppose that $\gcd(3,q-1)=3$. 
Then $3 \mid q^2+q+1$ so $q^2+q+1$ is not a prime. 
Thus the largest prime divisor is at most $\max\left\{\frac{q^2+q+1}{3}, q+1\right\} = \frac{q^2+q+1}{3}$, 
since $q>2$, whence $\frac{m}{2} = \frac{q^2}{2} > \frac{q^2+q+1}{3}$, again since $q>2$.

We thus see that the only possibility for $\psl_3(q)$ to have a large prime divisor is when $q^2+q+1$ is a prime number. 
When that is the case, $q^2+q+1 = m$ is in fact the only large prime.

Finally, we address the $n=2$ case. 
The arguments are similar, so we omit the details. 
The only possibilities for $\psl_2(q)$ to have at least one large prime divisor for its order is when $q$ is an odd prime, 
or $\frac{q+1}{2}$ is a prime but $q$ is not ($q$ odd), 
or when $q$ is even and $q+1$ is a Fermat prime. 
Further, it is only possible that $\psl_2(q)$ has a second large prime divisor for its order 
if $q=p$ is a prime such that $\frac{p+1}{2}$ is also a prime.
\subsection{Unitary case}
For unitary groups we need not consider the case $n=2$ as $\psu_2(q) \cong \psl_2(q)$. 
Suppose first that $n\geq 5$. 
We know that $m(\psu_n(q)) \geq q^{\lfloor\frac{n^2}{4}\rfloor}$ and we have already argued that 
$\lfloor \frac{n^2}{4}\rfloor \geq n+1$ for $n$ in that range. 
From the formula 
\begin{equation}
\lvert\psu_n(q)\rvert = \frac{1}{\gcd(n,q+1)}q^{n(n-1)/2}\prod_{i=2}^{n}(q^i-(-1)^i)\,,
\end{equation}
we see that the largest prime divisor of the order is at most $\frac{q^{n-1}+1}{q+1}$, if $n$ is even, 
and at most $\frac{q^{n}+1}{q+1}$ if $n$ is odd, thus at most $\frac{q^{n}+1}{q+1}$ in either case.
Therefore
\begin{equation}
\frac{m(\psu_n(q))}{2} \geq \frac{q^{n+1}}{2}  > \frac{q^n+1}{q+1}\,,
\end{equation}
and no large prime divisor occurs if $n \geq 5$. 
If $n=4$ then $m(\psu_4(q)) \geq q^4$, while the largest prime divisor of $\lvert\psu_4(q)\rvert$ is at most $q^2-q+1$. 
Clearly, $\frac{q^4}{2} > q^2-q+1$, thus we need only consider the case $n=3$.

For $n=3$ note that 
\begin{equation}
\lvert\psu_3(q)\rvert = \frac{1}{\gcd(3,q+1)}q^{3}(q-1)(q+1)^2(q^2-q+1)\,,    
\end{equation}
so the largest prime divisor is at most $q^2-q+1$. On the other hand,
$m(\psu_3(q)) \geq q^2$ and we may, of course, assume that $q>2$. 
Now if $q^2-q+1=p$ is a prime, then $\frac{m}{2}<p$, 
but the second largest prime divisor of the order then cannot exceed $q+1<\frac{q^2}{2} \leq \frac{m}{2}$.
If, however, $q^2-q+1$ is not a prime, then the largest prime divisor 
is at most $\max\left\{\frac{q^2-q+1}{3}, q+1\right\}$ (since $q^2-q+1$ is odd), 
which is easily seen to be less than $\frac{q^2}{2} \leq \frac{m}{2}$.

In conclusion, the only unitary groups with at least one (and thus exactly one) large prime divisor for their order are $\psu_3(q)$, 
where $q$ is such that $q^2-q+1$ is a prime number.
\subsection{Symplectic case}
Symplectic groups can be dealt with quickly. For one, the maximal abelian order of $\psp_{2n}(q)$ 
admits a lower bound which does not involve the floor function, thus resulting in less cases requiring consideration. 
As we are about to see, no symplectic group (which is not linear) has a large prime divisor for its order. 
For $n=1$ we have $\psp_2(q) = \psl_2(q)$, thus we restrict to $n \geq 2$. 
For $n=q=2$ we have the exceptional isomorphism $\psp_4(2) \cong S_6$ and $m(S_6)=9$, 
so $|S_6|$ has a single large prime divisor. Thus if $n=2$ we require $q>2$. Now recall that
\begin{equation}
\lvert\psp_{2n}(q)\rvert = \frac{1}{\gcd(2,q-1)}q^{n^2}\prod_{i=1}^{n}(q^{2i}-1)\,.
\end{equation}
It readily follows that the largest prime divisor of $\lvert \psp_{2n}(q) \rvert$ is at most $q^n+1$. 
Note that $q^n+1$ can, in fact, be a prime if $n$ is a power of 2 and $q^n+1$ is a Fermat prime. 
However, for all $n \geq 2$ we have $m(\psp_{2n}(q)) \geq q^{\frac{n(n+1)}{2}}$ 
with equality if $n \geq 3$ (for $n=2$ the inequality may be strict). It is now trivial to verify that 
\begin{equation}
m(\psp_{2n}(q)) \geq q^{\frac{n(n+1)}{2}} > q^n+1\,,
\end{equation}
since the only case where that inequality may fail, i.e. $n=2$, requires $q=2$.
\subsection{Orthogonal case}
Finally, we treat orthogonal groups and we address the odd-dimension case first. Recall that
\begin{equation}
\lvert O_{2n+1}(q)\rvert = \frac{1}{\gcd(2,q-1)}q^{n^2}\prod_{i=1}^{n}(q^{2i}-1)\,,
\end{equation}
where $n \geq 3$. 
Clearly, the size of the largest prime divisor of the order is at most $q^n+1$, 
whereas $m(O_{2n+1}(q)) \geq q^{\frac{n(n-1)}{2}+1}$ and it is easy to see that $m/2 > q^n+1$ always.

In even dimensions of plus-type we have
\begin{equation}
\lvert O_{2n}^{+}(q)\rvert = \frac{1}{\gcd(4,q^n-1)}q^{n(n-1)}(q^n-1)\prod_{i=1}^{n-1}(q^{2i}-1)\,,
\end{equation}
whence $q^{n}-1<q^n+1$ is again an upper bound for the largest prime divisor. 
Again, $m(O_{2n}^{+}(q)) \geq q^{\frac{n(n-1)}{2}+1}$ and it follows that $m/2 > q^n+1$ as before.

In even dimensions of minus-type we have
\begin{equation}
\lvert O_{2n}^{-}(q)\rvert = \frac{1}{\gcd(4,q^n+1)}q^{n(n-1)}(q^n+1)\prod_{i=1}^{n-1}(q^{2i}-1)\,.
\end{equation}
Once more $q^n+1$ is an upper bound for the largest prime, but this time 
\begin{equation}
m(O_{2n}^{-}(q)) \geq q^{\frac{(n-1)(n-2)}{2}+2} \, ;
\end{equation}
equality is possible if $q$ is odd. However, it is readily seen that 
\begin{equation}
    \frac{(n-1)(n-2)}{2}+2 \geq n+1
\end{equation}
for all $n \geq 4$, whence 
\begin{equation}
    m(O_{2n}^{-}(q)) \geq q^{\frac{(n-1)(n-2)}{2}+2} \geq q^{n+1} > q^n+1\,.
\end{equation}
Thus the orthogonal case contributes no additional groups.
\subsection{Exceptional groups}
Exceptional groups of Lie type are rather easy to deal with, 
since we can obtain exact estimates for the size of the largest prime divisor of the order 
(that is, the expressions involve only the size of the field).

Let us begin with $G_2(q)$. We have 
\begin{equation}
    |G_2(q)| = q^6(q^6-1)(q^2-1) = q^6(q-1)^2(q+1)^2(q^2-q+1)(q^2+q+1)\,,
\end{equation}
thus the largest prime divisor of the order is at most $q^2+q+1$. 
Now $m(G_2(q))$ is at least either $q^4$ or $q^3$ depending as the characteristic of the field is 3 or not. 
We can assume that $q>2$, as $G_2(2)$ is not simple; in fact, $G_2(2) \cong \Aut(\psu_3(3))$ is almost simple 
and $m(G_2(2)) =16$, while $|G_2(2)| = 2^6.3^3.7$. Then
\begin{equation}
    m(G_2(q)) \geq q^3 > 2\left(\frac{q^3-1}{q-1}\right) = 2(q^2+q+1)\,,
\end{equation}
for $q \geq 3$ and we see that groups of type $G_2(q)$ have no large prime divisors. 
Similarly for the small Ree groups $^2G_2(q)$, where $q$ is an odd power of 3, 
we see that the largest prime divisor is at most 
\begin{equation}
\frac{1}{2}(q^2-q+1) < \frac{1}{2}q^2 \leq \frac{m}{2}\,.    
\end{equation}

Next, for groups of type $F_4(q)$ we have
\begin{align*}
    |F_4(q)|    &= q^{24}(q^{12}-1)(q^8-1)(q^6-1)(q^2-1) \\
                &= q^{24}(q-1)^4  (q+1)^4 (q^2+1)^2 (q^2-q+1)^2 (q^2+q+1)^2 (q^4+1) (q^4-q^2+1)\,,
\end{align*}
thus the largest prime divisor of $|F_4(q)|$ is at most $q^4+1$. 
On the other hand, $m(F_4(q)) \geq q^9$ in all cases and we see that $\frac{m}{2} > q^4+1$ for all prime powers $q$.

As for the large Ree groups $^2F_4(q)$, where $q$ is an odd power of 2, their factorised order is
\begin{equation}
\left|^2F_4(q)\right| = q^{12} (q-1)^2 (q+1)^2 (q^2+1)^2 (q^2-q+1)(q^4-q^2+1)\,,
\end{equation}
and we see that $q^4-q^2+1$ is an upper bound for the size of the largest prime divisor of the order, 
while $m(^2F_4(q)) \geq q^5 > 2(q^4-q^2+1)$.

As for the groups $E_6(q)$, $E_7(q)$, and $E_8(q)$ the largest prime divisors for their orders are, respectively, 
at most $q^6+q^3+1$, $\frac{q^7-1}{q-1}$, and $q^8+q^7-q^5-q^4-q^3+q+1$. 
On the other hand, $m(E_6(q)) \geq q^{16}$, $m(E_7(q)) \geq q^{27}$, and $m(E_8(q)) \geq q^{36}$. 
For the groups $^2E_6(q)$ we find that $m(^2E_6(q)) \geq q^{12}$, while the largest prime divisor is at most $q^6-q^3+1$.
It is clear that these families contribute no additional groups with at least one large prime divisor.

Next, we consider the groups $^3D_4(q)$. Their (factorised) order is
\begin{equation}
\left|^3D_4(q)\right|   = q^{12}(q-1)^2 (q+1)^2 (q^2-q+1)^2(q^2+q+1)^2 (q^4-q^2+1)\,. 
\end{equation}
We see that the largest prime divisor is at most $q^4-q^2+1$, while 
\begin{equation}
m(^3D_4(q)) \geq q^5 > 2(q^4-q^2+1).    
\end{equation}

The only case left is the Suzuki groups $^2B_2(q)$, where $q=2^{2n+1}$ and $n \geq 1$. 
Notice that 
\begin{equation}
    \left|^2B_2(q)\right| = q^2(q^2 + 1) (q - 1) = q^2 (q - 1) (q + \sqrt{2q} + 1)(q - \sqrt{2q} + 1)\,,
\end{equation}
thus the largest prime divisor is at most $q + \sqrt{2q} + 1$.   If this number is not prime, 
then the largest prime divisor has order less than $q$.  

We claim that $m(^2B_2(q))=2q$.   Let $A$ be a maximal abelian subgroup.   Then
the centraliser of an involution is a Sylow $2$-subgroup and so either $A$ is a $2$-group
or has odd order.   Since every odd Sylow subgroup is cyclic, any odd order abelian
subgroup is cyclic and since we know all subgroups, we see that the maximal
order of an element is $q + \sqrt{2q} + 1$.   Since the centre of a Sylow $2$-subgroup
is elementary abelian of order $q$, there certainly are abelian subgroups of order $2q$.
We note that the centre of a Sylow $2$-subgroup has index $2$ in the centraliser
of any element of order $4$ (there is a unique conjugacy class of subgroups of order $4$
and this can be seen already in $\Sp_4(q)$--see \cite{LS}).   This proves the claim.  

Thus, we see that there exists a prime divisor $p > m/2$ if and only if
$q + \sqrt{2q} + 1$ is prime and in that case there is only $1$.

We are now able to summarise our findings as follows.  In the first result
we only list the non-abelian simple groups with one prime divisor.   There are
examples of almost simple groups as well (but the socle must have at least
one large prime divisor and indeed exactly one by the results below).  For
example if $p \geq 5$ is prime, then $m(\pgl_2(p)) =p+1$.


\begin{theorem}\label{thm:1prime}   Suppose that $G$ is a
non-abelian finite simple group with
exactly one prime divisor $p > m(G)/2$.  Then one of the following holds:
\begin{enumerate}
\item  $G \cong A_6$ or $A_7$;
\item $G$ is isomorphic to one of $\M_{11}$, $\M_{12}$, $\M_{22}$ or $\M_{23}$;
\item  $G \cong \psl_3(q)$ with $p=(q^3-1)/(q-1)$;
\item  $G \cong \psu_3(q), q > 2$ with $p=(q^3+1)/(q+1)$;
\item  $G \cong \psl_2(p)$, $p \geq 5$ with $(p+1)/2$ not prime; or
 \item  $G \cong {^2}B_2(q), q = 2^{2n+1} \geq 8$ with $p=q + \sqrt{2q} + 1$.
\end{enumerate}
\end{theorem}

\begin{theorem}\label{thm:C}
The only non-abelian finite simple groups with at least two, and 
thus exactly two, large prime divisors for their orders are 
$A_5$, $\J_1$, $\J_3$ and $\psl_2(p)$, where $p>5$ is a prime 
such that $\frac{p+1}{2}$ is also a prime.
\end{theorem}

At this point we need to address the ``pure" almost simple case and we do this next.

\begin{corollary}\label{cor:first}
Suppose that $S \leq G \leq \Aut(S)$, where $S$ is a non-abelian 
finite simple group and that $p, q$ are primes 
dividing $|G|$ with $q > p > m(G)/2$. Then $G \cong S$ and thus $G$ is one of the groups in Theorem~\ref{thm:C}.
\end{corollary}
\begin{proofof}
In the proof of Theorem~\ref{thm:B} we saw that $p, q$ must, in fact, divide $|S|$, since $S = \Soc(G)$. 
Thus $S$ is necessarily one of the groups in Theorem~\ref{thm:C}. Now $\J_1$ has no outer automorphisms, 
so if $S=\J_1$ then $G=\J_1$, while $S_5 \cong \Aut(A_5)$ and $m(S_5)=6$ but 5 is the only large prime dividing $|S_5|$. 
Thus if $S=A_5$ then $G=A_5$. As for $\J_3$, we have $\Aut(\J_3) \cong \J_3.2$ and we find--using {\GAP} again--that 
$m(\Aut(\J_3))=34$ whence $\Aut(\J_3)$ has only one large prime divisor for its order. Thus if $S=\J_3$ then $G=\J_3$.

Therefore, it suffices to argue the case $S=\psl_2(p)$ with $p>5 $ a prime such that $\frac{p+1}{2}$ is also a prime. 
Since $p$ is prime, we have $\Aut(\psl_2(p)) \cong \pgl_2(p)$, while $m(\pgl_2(p)) = p+1$ whence $p$ is the unique 
prime divisor of $\lvert\pgl_2(p)\rvert$ which is strictly larger than the prime $\frac{m}{2}= \frac{p+1}{2}$. 
In that case too, therefore, it follows that if $S = \psl_2(p)$ then $G = \psl_2(p)$.
\end{proofof}

We are now ready to argue the case of equality in Theorem~\ref{thm:A}. 
Recall that $g(n)$ is defined as the product of all prime powers at most $n$.

\begin{corollary}\label{cor:second}
Let $G$ be a non-trivial group such that $|G|=g(m)$, where $m=m(G)$ is the maximal abelian order of $G$. 
Then $m \in \{2,3,4,6\}$ and $G$ is isomorphic to $S_2$, $S_3$, $S_4$, and $S_5$ respectively.
\end{corollary}
\begin{proofof}
Nagura \cite{nagura} has shown that for all $x \geq 25$ there is at least on prime in the interval $(x, 6/5x)$. 
By a twofold application of this result, we ensure the existence of at least two primes in the range 
$(x,2x)$ for all $x \geq 25$. 
By inspecting small values of $m$, we may therefore deduce that there are at least two primes in the interval 
$(m/2,m]$ for all $m \in \mathbb{N}_{\geq 3}$, except when $m \in \{4,6,10\}$.

It is clear that when $m=2$ or $m=3$ the only groups of order $g(2)$, $g(3)$ respectively, are $S_2$, $S_3$. 
When $m=4$ or $m=6$, we employ {\GAP} to argue that the only groups of size $g(4)$, $g(6)$ are 
$S_4$ and $S_5$, respectively. Lastly, when $m=10$ we argue that a group of order $g(10)$ has an abelian 
subgroup of size $>10$. Note that if $|G|=g(10)$, then a Sylow $2$-subgroup of $G$ has order $2^6$. 
Using {\GAP}, we find that all groups of order $2^6$ have an abelian subgroup of size at least $2^4$.

By Corollary~\ref{cor:first}, the proof will be complete provided that we can show that each of the groups  
$A_5$, $\J_1$, $\J_3$ and $\psl_2(p)$, where $p>5$ is a prime such that $\frac{p+1}{2}$ is also a prime, 
does not have order $g(m)$. For  $A_5$, $\J_1$, and $\J_3$ this is clear by a mere inspection 
and it is just as easy to see that 
\begin{equation}
\lvert \psl_2(p)\rvert = \frac{1}{2}p(p-1)(p+1) \neq g(p)\,,   
\end{equation}
which is equivalent to 
\begin{equation}
	(p-1)\left(\frac{p+1}{2}\right) \neq g(p-1)\,.    
\end{equation}
Since $p \geq 13$, there are at least two primes in the range  $(\frac{p-1}{2},p-1]$ and one of those is $\frac{p+1}{2}$. 
A second large prime divides $g(p-1)$ but not $p-1$, and thus our proof is complete.
\end{proofof}

As we promised in the Introduction, we are now able to strengthen Theorem~\ref{thm:A}.
For the purposes of the next theorem only, define $h(n)$ to be the product of all primes 
at most $n$ which are strictly larger than $n/2$.

\begin{theorem}\label{thm:g/h}
Let $G$ be a finite group not isomorphic to one of the groups in Theorem~\ref{thm:C} (having two large prime divisors), 
and suppose that $m=m(G)$. Then there is a prime $p \in (m/2,m]$ such that $|G|$ divides $p \cdot g(m)/h(m)$. 
As a consequence, for all finite groups $G$ we have
\begin{equation}
	|G| \leq m  \frac{g(m)}{h(m)}\,,
\end{equation}
except when $G \cong S_3$, or $G \cong A_5$.
\end{theorem}
The proof is only a matter of case-checking, so we omit the details.
\section{Asymptotic estimate for the upper bound}
Fix a positive integer $n$ and observe that the function $g(n)$, which is the product of all prime powers at most $n$, can be written as
\begin{equation}
g(n) = \prod_{p\leq n} p^{\frac{\xi_{p}(1+\xi_{p})}{2}}\,,
\end{equation}
where $\xi_{p}$ is the unique positive integer such that 
\begin{equation}
p^{\xi_{p}}\leq n <p^{1+\xi_{p}}\,.
\end{equation}
From this we note that 
\begin{equation}\label{eq:parena}
\xi_{p} \leq \frac{\log n}{\log p}\,,
\end{equation}
and also that
\begin{equation}\label{eq:parduo}
\sqrt{n}<p\leq n
\Rightarrow \xi_{p} = 1.
\end{equation}
Recall that, by definition,
\begin{equation}
h(n) = \prod_{\frac{n}{2} < p \leq n} p\, .
\end{equation}

Our goal in this section is to estimate the order of magnitude of the upper bound in Theorem~\ref{thm:g/h} as follows.
\begin{lemma}
Let $f(n) \coloneqq n \frac{g(n)}{h(n)}$, where $g(n)$ and $h(n)$ are defined as above. Then
\begin{equation}\label{eq:willsho}
f(n)=e^{\frac{n}{2}+O\left(\frac{n}{\log n}\right)}\,.
\end{equation}
In particular, for every $\varepsilon \in (0,1)$ there exists $n_0=n_0(\varepsilon)$ such that 
if $n\geq n_0$ then 
\begin{equation}
(e-\varepsilon)^{\frac{n}{2}} < f(n) < (e+\varepsilon)^{\frac{n}{2}}\,.
\end{equation}
\end{lemma}
\begin{proofof}
We begin by noting that
\begin{equation}\log g(n)=
\sum_{p\leq n} {\frac{\xi_{p}(1+\xi_{p})}{2}} \log p\,,
\end{equation} 
hence by~\eqref{eq:parduo} we have 
\begin{align*}
\log g(n) 	&= \sum_{p\leq \sqrt{n}} {\frac{\xi_{p}(1+\xi_{p})}{2}} \log p + \sum_{\sqrt{n} < p\leq n} \log p \\
		&= \sum_{p\leq \sqrt{n}} {\frac{\xi_{p}(1+\xi_{p})}{2}} \log p + \sum_{\sqrt{n} < p\leq \frac{n}{2}} \log p + \sum_{\frac{n}{2} < p\leq n} \log p\,.
\end{align*}
Therefore
\begin{equation}
\log g(n) - \log h(n) 	= \sum_{p\leq \sqrt{n}} {\frac{\xi_{p}(1+\xi_{p})}{2}} \log p + \sum_{\sqrt{n} < p\leq \frac{n}{2}} \log p\,.
\end{equation}

A form of the Prime Number Theorem is
\begin{equation}
\sum_{p\leq n} \log p=n+O\left(\frac{n}{\log n}\right)\,,
\end{equation}
whence 
\begin{align*}
\sum_{\sqrt{n}<p\leq \frac{n}{2}} \log p = \sum_{p\leq \frac{n}{2}} \log p - \sum_{p\leq \sqrt{n}} \log p 	&= \frac{n}{2}+O\left(\frac{n}{\log n}\right)+O(\log n \sqrt{n})	\\
																				&= \frac{n}{2}+O\left(\frac{n}{\log n}\right)\,.
\end{align*}
This shows that 
\begin{equation}\label{eq:steff}
\log g(n) - \log h(n) = \sum_{p\leq \sqrt{n}} {\frac{\xi_{p}(1+\xi_{p})}{2}} \log p + \frac{n}{2} + O\left(\frac{n}{\log n}\right)\,.
\end{equation}
Now we prove that the last sum over $p$ is much smaller than $\frac{n}{2}$.
Indeed, by~\eqref{eq:parena}
we obtain
\begin{align*}
\sum_{p\leq \sqrt{n}} \xi_{p}{\frac{(1+\xi_{p})}{2}} \log p  	&\leq  \sum_{p\leq \sqrt{n}} \frac{\log n}{\log p}  \left(\frac{1}{2}+\frac{\log n}{2 \log p}\right) \log p					\\
												&\leq  \sum_{p\leq \sqrt{n}} \frac{\log n}{\log p}  \left(\frac{1}{2}\frac{\log n}{\log p}+\frac{\log n}{2 \log p}\right) \log p	\\
												&\leq  (\log n)^2 \sum_{p\leq \sqrt{n}}   \frac{1}{\log p}													\\
												&\leq (\log n)^2 \sqrt{n}\,.
\end{align*}
Injecting this into~\eqref{eq:steff} yields 
\begin{equation}
\log g(n) - \log h(n) = \frac{n}{2} + O\left(\frac{n}{\log n}\right),
\end{equation}
whence
\begin{equation}
\log f(n) = \log n + \log g(n) - \log h(n) = \frac{n}{2} + O\left(\frac{n}{\log n}\right),
\end{equation}
completing the proof of~\eqref{eq:willsho}.
\end{proofof}

\section{Afterword}
Upon completing the writing of this note, it came to our attention that Erd\H{o}s and Straus~\cite{erdosstraus}
had also considered the question what size we can guarantee for an abelian subgroup of a group in terms of
the group's order. In fact, their Theorem 1.1, although less precise in its conclusion, is a close neighbour of
our own Theorem~\ref{thm:A}. Laci Pyber \cite{pyber} then picked up the problem and essentially resolved it
(Erd\H{o}s put a copy of his paper in Pyber's mailbox). The conclusion of Pyber's main theorem is that there
is an absolute constant $c$ such that a finite group of order $n$ has an abelian subgroup of order
$2^{c \sqrt{\log n}}$. By a result of Ol'shanski\u{i}~\cite{olshanskii}, this estimate is best possible up 
to the constant~$c$. Of course, Pyber's theorem is asymptotically much better than Theorem~\ref{thm:g/h}.
It also has as a consequence that there can only exist finitely many groups $G$ and $m = m(G) \in \mathbb{N}$
such that $|G| = g(m)$, as per Corollary~\ref{cor:second}, although it does not give more
precise information about $m(G)$ or $G$.
\appendix
\section{GAP code}
Here we present the snippet of code we have used to determine (when possible) 
the largest abelian order of a group, along with a few sample calls.
\begin{verbatim}
gap> LoadPackage("AtlasRep");
true
gap> LargestAbelianOrder:=function(g)
> local list,cls,abelcls;
> cls:=ConjugacyClassesSubgroups(g);
> abelcls:=Filtered(cls, c -> IsAbelian(Representative(c)));
> list:=List(abelcls, c -> Order(Representative(c)));
> return Maximum(Set(list));
> end;;
gap> J1:=AtlasGroup("J1");
<permutation group of size 175560 with 2 generators>
gap> LargestAbelianOrder(J1);        
19
gap> J3:=AtlasGroup("J3");   
<permutation group of size 50232960 with 2 generators>
gap> LargestAbelianOrder(J3);                                  
27
\end{verbatim}

\bibliographystyle{amsalpha}

\end{document}